\documentclass[12pt, reqno]{amsart}
\usepackage{amsmath, amsthm, amscd, amsfonts, amssymb, graphicx, color}
\usepackage[bookmarksnumbered, colorlinks, plainpages]{hyperref}

\hypersetup{colorlinks=true,linkcolor=red, anchorcolor=green, citecolor=cyan, urlcolor=red, filecolor=magenta, pdftoolbar=true}
\newtheorem{theorem}{Theorem}[section]

\theoremstyle{definition}

\theoremstyle{remark}

\numberwithin{equation}{section}
\input{tcilatex}

\begin{document}
\title[Existence of solutions for a class of semilinear system]{\textbf{On
the} \textbf{radial solutions of a system with weights under the
Keller--Osserman condition}}
\author[D.-P. Covei]{Dragos-Patru Covei$^{1}$ }
\address{$^{1}$Department of Applied Mathematics, The Bucharest University
of Economic Studies\\
Piata Romana, 1st district, postal code: 010374, postal office: 22, Romania}
\email{\textcolor[rgb]{0.00,0.00,0.84}{coveid@yahoo.com}}
\date{Received: xxxxxx; Revised: yyyyyy; Accepted: zzzzzz. \\
\indent $^{*}$Corresponding author}

\begin{abstract}
In this paper, we establish conditions on the functions $p_{1}$ and $p_{2}$
that are necessary and sufficient for the existence of positive solutions,
bounded and unbounded, of the semilinear elliptic system%
\begin{equation*}
\left\{ 
\begin{array}{l}
\Delta u=p_{1}\left( \left\vert x\right\vert \right) f_{1}\left( u,v\right) 
\text{ for }x\in \mathbb{R}^{N}\text{ (}N\geq 3\text{), } \\ 
\Delta v=p_{2}\left( \left\vert x\right\vert \right) f_{2}\left( u,v\right) 
\text{ for }x\in \mathbb{R}^{N}\text{ (}N\geq 3\text{),}%
\end{array}%
\right. 
\end{equation*}%
where $p_{1}$, $p_{2}$, $f_{1}$ and $f_{2}$ are continuous functions
satisfying certain properties.
\end{abstract}

\maketitle

\setcounter{page}{1}


\let\thefootnote\relax\footnote{%
Copyright 2016 by the Journal.}

\section{Introduction}

The study of existence of large solutions for semilinear elliptic systems of
the form%
\begin{equation}
\left\{ 
\begin{array}{l}
\Delta u=p_{1}\left( \left\vert x\right\vert \right) f_{1}\left( u,v\right) 
\text{ for }x\in \mathbb{R}^{N}\text{ (}N\geq 3\text{), } \\ 
\Delta v=p_{2}\left( \left\vert x\right\vert \right) f_{2}\left( u,v\right) 
\text{ for }x\in \mathbb{R}^{N}\text{ (}N\geq 3\text{),}%
\end{array}%
\right.  \label{11}
\end{equation}%
goes back to the pioneering papers by Keller \cite{K} and Osserman \cite{O}.
\ From the results of \cite[Lemma 3, pp. 1643]{O} we know that, for a given
positive, continuous and nondecreasing function $f$, the semilinear elliptic
partial differential inequality 
\begin{equation}
\Delta u\geq f\left( u\right) \text{ in }\mathbb{R}^{N},  \label{iko}
\end{equation}%
possesses a large solution $u$ if and only if the nowadays called
Keller-Osserman condition holds, i. e. 
\begin{equation}
KO_{f}:=\int_{1}^{\infty }\frac{1}{\sqrt{F\left( s\right) }}ds=+\infty \text{%
, (}F\left( s\right) =\int_{0}^{s}f\left( t,t\right) dt\text{).}  \label{ko}
\end{equation}%
Here, we extend the existence result to the case of systems where the
functions $p_{1}$ and $p_{2}$ are spherically symmetric. More generally,
however, we are interested in the influence of the functions $p_{1}$ and $%
p_{2}$ on existence results. Here, instead of fixing the conditions on $%
p_{1} $ and $p_{2}$, we fix the conditions on $f_{1}$ and $f_{2}$ and
determine sufficient condition on $p_{1}$ and $p_{2}$ that ensure that (\ref%
{11}) has an entire solution and whether such solutions are bounded or
unbounded and, perhaps, large. Finally, we note that the study of large or
bounded solutions for (\ref{11}) when $KO_{f}$ $<\infty $ or $KO_{f}=\infty $
has been the subject of many articles. See, for example, the author \cite%
{CD2}, Lair \cite{LAIR}, Nehari \cite{NE}, Rhee \cite{RH} and Redheffer \cite%
{RR}, and their references.

The problems of the form (\ref{11}) and (\ref{iko}) are drawn by the
mathematical modelling of many natural phenomena related to steady-state
reaction-diffusion, subsonic fluid flows, electrostatic potential in a shiny
metallic body inside or subsonic motion of a gas, automorphic functions
theory, geometry and control theory (see, for example, L. Bieberbach \cite%
{LB}, Grosse-Martin \cite{GR}, Diaz \cite{DZ}, Keller \cite{KG}, Lasry and
Lions \cite{LAS}, Matero \cite{JM}, Pohozaev \cite{PC4}\textit{,} Rademacher%
\textit{\ }\cite{R} and Smooke \cite{S}\textit{\ }for a more detailed
discussion). For example, reading the work of Lasry and Lions \cite{LAS}, we
can observe that such problems arise in stochastic control theory. The
controls are to be designed so that the state of the system is constrained
to some region. Finding optimal controls is then shown to be equivalent to
finding large solutions for a second order nonlinear elliptic partial
differential equation.

\section{The main results}

Let $a_{1},a_{2}\in \left( 0,\infty \right) $. We assume:

(P1)\quad $p_{1},p_{2}:\left[ 0,\infty \right) \rightarrow \left[ 0,\infty
\right) $ are spherically symmetric continuous functions (i.e.,\textit{\ }for%
\textit{\ }$r=\left\vert x\right\vert $\textit{\ we have }$p_{1}\left(
x\right) =p_{1}\left( r\right) $ and $p_{2}\left( x\right) =p_{2}\left(
r\right) $)\textit{;}

(C1)\quad $f_{1},f_{2}:\left[ 0,\infty \right) \times \left[ 0,\infty
\right) \rightarrow \left[ 0,\infty \right) $ are continuous, increasing in
each variables, $f_{1}\left( 0,0\right) f_{2}\left( 0,0\right) \geq 0$ and $%
f_{1}\left( s_{1},s_{2}\right) f_{2}\left( s_{1},s_{2}\right) >0$ whenever $%
s_{1},s_{2}>0$;

(C2)\quad there exist the continuous and increasing functions $\overline{%
f_{1}},$ $\overline{f_{2}}:\left[ 0,\infty \right) \rightarrow \left[
0,\infty \right) $ and $\overline{c}_{1},\overline{c}_{2}\in \left( 0,\infty
\right) $ such that 
\begin{eqnarray}
f_{1}\left( t_{1},t_{1}s_{1}\right) &\leq &\overline{c}_{1}f_{1}\left(
t_{1},t_{1}\right) \overline{f_{1}}\left( s_{1}\right) \text{, }\forall 
\text{ }t_{1},s_{1}>0,  \label{c22} \\
f_{2}\left( t_{2},t_{2}s_{2}\right) &\leq &\overline{c}_{2}f_{2}\left(
t_{2},t_{2}\right) \overline{f_{2}}\left( s_{2}\right) \text{, }\forall 
\text{ }t_{2},s_{2}>0.  \label{c222}
\end{eqnarray}

A simple example of $f_{1}$ and $f_{2}$ satisfying (C1) and (C2) is given by
\ $f_{1}\left( u,v\right) =v^{\alpha }$ and $f_{2}\left( u,v\right)
=u^{\beta }$ with $\alpha $, $\beta \in \left( 0,\infty \right) $. The
assumption (C2) is further discussed in the famous book of Krasnosel'skii
and Rutickii \cite{KR} (see also Rao and Ren \cite{RAO}).

Next, we introduce the following notations%
\begin{eqnarray*}
Z\left( r\right) &=&\int_{a_{1}+a_{2}}^{r}\frac{1}{f_{1}\left( t,t\right)
+f_{2}\left( t,t\right) }dt\text{, } \\
\mathcal{KO}_{f_{1}}\left( r\right) &=&\int_{a_{1}}^{r}\frac{1}{\sqrt{%
\int_{0}^{s}f_{1}\left( t,t\right) dt}}ds\text{, }\mathcal{KO}_{f_{1}}\left(
\infty \right) =\lim_{s\rightarrow \infty }\mathcal{KO}_{f_{1}}\left(
s\right) \text{, } \\
\mathcal{KO}_{f_{2}}\left( r\right) &=&\int_{a_{2}}^{r}\frac{1}{\sqrt{%
\int_{0}^{s}f_{2}\left( t,t\right) dt}}ds\text{, }\mathcal{KO}_{f_{2}}\left(
\infty \right) =\lim_{s\rightarrow \infty }\mathcal{KO}_{f_{2}}\left(
s\right) \text{, } \\
M_{i} &\geq &\max \left\{ 1,\frac{1}{a_{i}}\right\} \text{, }i=1,2\text{, }%
\varepsilon \in \left( 0,\infty \right) \text{, }\phi _{i}\left( s\right) =%
\underset{0\leq t\leq s}{\max }p_{i}\left( t\right) \text{,} \\
P_{i}\left( r\right) &=&\int_{0}^{r}z^{1-N}\int_{0}^{r}t^{N-1}p_{i}\left(
t\right) dtdz\text{, } \\
\overline{P}_{i}\left( r\right) &=&\sqrt{\overline{f_{i}}\left( M_{i}\left(
1+Z^{-1}\left( P_{1}\left( r\right) +P_{2}\left( r\right) \right) \right)
)\right) }\int_{0}^{r}\sqrt{\phi _{i}\left( s\right) }ds\text{, }\overline{P}%
_{i}\left( \infty \right) =\lim_{r\rightarrow \infty }\overline{P}_{i}\left(
r\right) \text{, } \\
\overline{P}_{i\varepsilon }\left( r\right) &=&\overline{f_{i}}\left(
M_{i}\left( 1+Z^{-1}\left( P_{1}\left( r\right) +P_{2}\left( r\right)
\right) \right) )\right) \int_{R}^{r}t^{1+\varepsilon }p_{i}\left( t\right)
dt\text{, }\overline{P}_{i\varepsilon }\left( \infty \right)
=\lim_{r\rightarrow \infty }\overline{P}_{i\varepsilon }\left( r\right) 
\text{,} \\
\underline{P}\left( r\right)
&=&\int_{0}^{r}y^{1-N}\int_{0}^{y}t^{N-1}p_{1}\left( t\right) f_{1}\left(
a_{1},a_{2}+f_{2}\left( a_{1},a_{2}\right) P_{2}\left( t\right) \right) dtdy%
\text{, }\underline{P}\left( \infty \right) =\lim_{r\rightarrow \infty }%
\underline{P}\left( r\right) \text{,} \\
\underline{Q}\left( r\right)
&=&\int_{0}^{r}y^{1-N}\int_{0}^{y}t^{N-1}p_{2}\left( t\right) f_{2}\left(
a_{1},a_{1}+f_{1}\left( a_{1},a_{2}\right) P_{1}\left( t\right) \right) dtdy%
\text{, }\underline{Q}\left( \infty \right) =\lim_{r\rightarrow \infty }%
\underline{Q}\left( r\right) ,
\end{eqnarray*}

The main results are:

\begin{theorem}
\label{th1}The system (\ref{11}) has one positive radial solution $\left(
u,v\right) \in C^{2}\left( \left[ 0,\infty \right) \right) \times
C^{2}\left( \left[ 0,\infty \right) \right) $ given \textit{that }$\mathcal{%
KO}_{f_{1}}\left( \infty \right) =\mathcal{KO}_{f_{2}}\left( \infty \right)
=\infty $ and \textrm{(P1)}, \textrm{(C1)}, \textrm{(C2) }hold true.
Moreover, if $\underline{P}\left( \infty \right) =\infty $ and $\underline{Q}%
\left( \infty \right) =\infty $ then 
\begin{equation*}
\lim_{\left\vert x\right\vert \rightarrow \infty }u\left( \left\vert
x\right\vert \right) =\infty \text{ and }\lim_{\left\vert x\right\vert
\rightarrow \infty }v\left( \left\vert x\right\vert \right) =\infty ,
\end{equation*}%
i.e. $\left( u,v\right) $ is an entire large solution.
\end{theorem}

\begin{theorem}
\label{th2}The system (\ref{11}) has one positive radial solution $\left(
u,v\right) \in C^{2}\left( \left[ 0,\infty \right) \right) \times
C^{2}\left( \left[ 0,\infty \right) \right) $ given that\textit{\ }$\mathcal{%
KO}_{f_{1}}\left( \infty \right) =\mathcal{KO}_{f_{2}}\left( \infty \right)
=\infty $ and \textrm{(P1)}, \textrm{(C1)}, \textrm{(C2)} hold true.
Moreover, if $r^{2N-2}p_{1}\left( r\right) $\textit{, }$r^{2N-2}p_{2}\left(
r\right) $\textit{\ are nondecreasing for large }$r$, $\overline{P}%
_{1\varepsilon }\left( \infty \right) <\infty $ and $\overline{P}%
_{2\varepsilon }\left( \infty \right) <\infty $ then 
\begin{equation*}
\lim_{\left\vert x\right\vert \rightarrow \infty }u\left( \left\vert
x\right\vert \right) <\infty \text{ and }\lim_{\left\vert x\right\vert
\rightarrow \infty }v\left( \left\vert x\right\vert \right) <\infty ,
\end{equation*}%
i.e. $\left( u,v\right) $ is an entire bounded solution. If, on the other
hand, (\ref{11}) has a nonnegative entire large solution, $%
r^{2N-2}p_{1}\left( r\right) $ respectively $r^{2N-2}p_{2}\left( r\right) $
are \textit{nondecreasing for large }$r$, then $p_{1}$ and $p_{2}$ satisfy 
\begin{equation*}
\overline{P}_{1\varepsilon }\left( \infty \right) =\overline{P}%
_{2\varepsilon }\left( \infty \right) =\infty .
\end{equation*}
\end{theorem}

\begin{theorem}
\label{th3}The system (\ref{11}) has one positive radial solution $\left(
u,v\right) \in C^{2}\left( \left[ 0,\infty \right) \right) \times
C^{2}\left( \left[ 0,\infty \right) \right) $ given that $\mathcal{KO}%
_{f_{1}}\left( \infty \right) =\mathcal{KO}_{f_{2}}\left( \infty \right)
=\infty $ and \textrm{(P1)}, \textrm{(C1)}, \textrm{(C2)} hold true.
Moreover, the following hold:

1)\quad If $r^{2N-2}p_{1}\left( r\right) $\textit{\ }is\textit{\
nondecreasing for large }$r$, $\overline{P}_{1\varepsilon }\left( \infty
\right) <\infty $ and $\underline{Q}\left( \infty \right) =\infty $ then 
\begin{equation*}
\lim_{\left\vert x\right\vert \rightarrow \infty }u\left( \left\vert
x\right\vert \right) <\infty \text{ and }\lim_{\left\vert x\right\vert
\rightarrow \infty }v\left( \left\vert x\right\vert \right) =\infty .
\end{equation*}

2)\quad If $r^{2N-2}p_{2}\left( r\right) $\textit{\ is nondecreasing for
large }$r$, $\underline{P}\left( \infty \right) =\infty $ and $\overline{P}%
_{2\varepsilon }\left( \infty \right) <\infty $ then 
\begin{equation*}
\lim_{\left\vert x\right\vert \rightarrow \infty }u\left( \left\vert
x\right\vert \right) =\infty \text{ and }\lim_{\left\vert x\right\vert
\rightarrow \infty }v\left( \left\vert x\right\vert \right) <\infty .
\end{equation*}
\end{theorem}

\begin{theorem}
\label{th4}The system (\ref{11}) has one positive bounded radial solution $%
\left( u,v\right) \in C^{2}\left( \left[ 0,\infty \right) \right) \times
C^{2}\left( \left[ 0,\infty \right) \right) $ given that $\overline{P}%
_{1}\left( \infty \right) <\mathcal{KO}_{f_{1}}\left( \infty \right) <\infty 
$, $\overline{P}_{2}\left( \infty \right) <\mathcal{KO}_{f_{2}}\left( \infty
\right) <\infty $, \textrm{(P1), (C1),\ (C2) }hold true. \ Moreover,%
\begin{equation*}
\left\{ 
\begin{array}{l}
a_{1}+\underline{P}\left( r\right) \leq u\left( r\right) \leq \mathcal{KO}%
_{f_{1}}^{-1}\left( \sqrt{2\overline{c}_{1}}\cdot \overline{P}_{1}\left(
r\right) \right) , \\ 
a_{2}+\underline{Q}\left( r\right) \leq v\left( r\right) \leq \mathcal{KO}%
_{f_{2}}^{-1}\left( \sqrt{2\overline{c}_{2}}\cdot \overline{P}_{2}\left(
r\right) \right) .%
\end{array}%
\right. \text{ }
\end{equation*}
\end{theorem}

\begin{theorem}
\label{th5}Assume \textrm{(P1), (C1) }and\textrm{\ (C2) }hold true. The
following hold true:

i)\quad The system (\ref{11}) has one positive radial solution $\left(
u,v\right) \in C^{2}\left( \left[ 0,\infty \right) \right) \times
C^{2}\left( \left[ 0,\infty \right) \right) $ such that $\lim_{r\rightarrow
\infty }u\left( r\right) =\infty $ and $\lim_{r\rightarrow \infty }v\left(
r\right) <\infty $ given that $\mathcal{KO}_{f_{1}}\left( \infty \right)
=\infty $, $\underline{P}\left( \infty \right) =\infty $ and $\overline{P}%
_{2}\left( \infty \right) <\mathcal{KO}_{f_{2}}\left( \infty \right) <\infty 
$.

ii)\quad The system (\ref{11}) has one positive radial solution $\left(
u,v\right) \in C^{2}\left( \left[ 0,\infty \right) \right) \times
C^{2}\left( \left[ 0,\infty \right) \right) $ such that $\lim_{r\rightarrow
\infty }u\left( r\right) <\infty $ and $\lim_{r\rightarrow \infty }v\left(
r\right) =\infty $ given that $\overline{P}_{1}\left( \infty \right) <%
\mathcal{KO}_{f_{1}}\left( \infty \right) <\infty $ and $\mathcal{KO}%
_{f_{2}}\left( \infty \right) =\infty $, $\underline{Q}\left( \infty \right)
=\infty $.
\end{theorem}

\section{Proofs of the main results}

Radial solutions of the system (\ref{11}) solve 
\begin{equation}
\left\{ 
\begin{array}{l}
\left( r^{N-1}u^{\prime }\left( r\right) \right) ^{\prime
}=r^{N-1}p_{1}\left( r\right) f_{1}\left( u\left( r\right) ,v\left( r\right)
\right) \text{ for }0\leq r<\infty , \\ 
\left( r^{N-1}v^{\prime }\left( r\right) \right) ^{\prime
}=r^{N-1}p_{2}\left( r\right) f_{2}\left( u\left( r\right) ,v\left( r\right)
\right) \text{ for }0\leq r<\infty .%
\end{array}%
\right.  \label{ss1}
\end{equation}%
We rewrite the system (\ref{ss1}) as 
\begin{equation}
\left\{ 
\begin{array}{l}
u\left( r\right) =u\left( 0\right)
+\int_{0}^{r}t^{1-N}\int_{0}^{t}s^{N-1}p_{1}\left( s\right) f_{1}\left(
u\left( s\right) ,v\left( s\right) \right) dsdt, \\ 
v\left( r\right) =v\left( 0\right)
+\int_{0}^{r}t^{1-N}\int_{0}^{t}s^{N-1}p_{2}\left( s\right) f_{2}\left(
u\left( s\right) ,v\left( s\right) \right) dsdt.%
\end{array}%
\right.  \label{ss}
\end{equation}%
We define the sequences $\left\{ u_{n}\right\} _{n\geq 0}$ and $\left\{
v_{n}\right\} _{n\geq 0}$ on $\left[ 0,\infty \right) $ iteratively by: 
\begin{equation}
\left\{ 
\begin{array}{l}
u_{0}=u\left( 0\right) =a_{1},v_{0}=v\left( 0\right) =a_{2}\text{,} \\ 
u_{n}\left( r\right)
=a_{1}+\int_{0}^{r}t^{1-N}\int_{0}^{t}s^{N-1}p_{1}\left( s\right)
f_{1}\left( u_{n-1}\left( s\right) ,v_{n-1}\left( s\right) \right) dsdt\text{
for }r\geq 0, \\ 
v_{n}\left( r\right)
=a_{2}+\int_{0}^{r}t^{1-N}\int_{0}^{t}s^{N-1}p_{2}\left( s\right)
f_{2}\left( u_{n-1}\left( s\right) ,v_{n-1}\left( s\right) \right) dsdt\text{
for }r\geq 0.%
\end{array}%
\right.  \label{recs}
\end{equation}%
We show that $\left\{ u_{n}\right\} _{n\geq 0}$ and $\left\{ v_{n}\right\}
_{n\geq 0}$ are nondecreasing on $\left[ 0,\infty \right) $. To see this,
express 
\begin{eqnarray*}
u_{1}\left( r\right)
&=&a_{1}+\int_{0}^{r}t^{1-N}\int_{0}^{t}s^{N-1}p_{1}\left( s\right)
f_{1}\left( u_{0}\left( s\right) ,v_{0}\left( s\right) \right) dsdt \\
&=&a_{1}+\int_{0}^{r}t^{1-N}\int_{0}^{t}s^{N-1}p_{1}\left( s\right)
f_{1}\left( a_{1},a_{2}\right) dsdt \\
&\leq &a_{1}+\int_{0}^{r}t^{1-N}\int_{0}^{t}s^{N-1}p_{1}\left( s\right)
f_{1}\left( u_{1}\left( s\right) ,v_{1}\left( s\right) \right)
dsdt=u_{2}\left( r\right) .
\end{eqnarray*}%
This proves that $u_{1}\left( r\right) \leq u_{2}\left( r\right) $. In the
same way $v_{1}\left( r\right) \leq v_{2}\left( r\right) $. By an induction
argument we get 
\begin{equation*}
u_{n}\left( r\right) \leq u_{n+1}\left( r\right) \text{ for any }n\in 
\mathbb{N}\text{ and }r\in \left[ 0,\infty \right)
\end{equation*}%
and%
\begin{equation*}
v_{n}\left( r\right) \leq v_{n+1}\left( r\right) \text{ for any }n\in 
\mathbb{N}\text{ and }r\in \left[ 0,\infty \right) .
\end{equation*}%
We show that the non-decreasing sequences $\left\{ u_{n}\right\} _{n\geq 0}$
and $\left\{ v_{n}\right\} _{n\geq 0}$ are bounded above on any compact
interval. By the monotonicity of $\left\{ u_{n}\right\} _{n\geq 0}$ and $%
\left\{ v_{n}\right\} _{n\geq 0}$ one get%
\begin{eqnarray}
\left[ r^{N-1}\left( u_{n}\left( r\right) \right) ^{\prime }\right] ^{\prime
} &=&r^{N-1}p_{1}\left( r\right) f_{1}\left( u_{n-1}\left( r\right)
,v_{n-1}\left( r\right) \right) \leq r^{N-1}p_{1}\left( r\right) f_{1}\left(
u_{n}\left( r\right) ,v_{n}\left( r\right) \right) ,  \label{gen1} \\
\left[ r^{N-1}\left( v_{n}\left( r\right) \right) ^{\prime }\right] ^{\prime
} &\leq &r^{N-1}p_{2}\left( r\right) f_{2}\left( u_{n}\left( r\right)
,v_{n}\left( r\right) \right) .  \label{gen2}
\end{eqnarray}%
By adding (\ref{gen1}) and (\ref{gen2}), one obtains%
\begin{eqnarray*}
\left[ r^{N-1}\left( u_{n}\left( r\right) +v_{n}\left( r\right) \right)
^{\prime }\right] ^{\prime } &\leq &r^{N-1}p_{1}\left( r\right) f_{1}\left(
u_{n}\left( r\right) ,v_{n}\left( r\right) \right) +r^{N-1}p_{2}\left(
r\right) f_{2}\left( u_{n}\left( r\right) ,v_{n}\left( r\right) \right) \\
&\leq &r^{N-1}\left( p_{1}\left( r\right) +p_{2}\left( r\right) \right)
\left( \left( f_{1}+f_{2}\right) \left( u_{n}\left( r\right) +v_{n}\left(
r\right) ,u_{n}\left( r\right) +v_{n}\left( r\right) \right) \right) .
\end{eqnarray*}%
Integration leads to%
\begin{equation}
\frac{\left( u_{n}\left( r\right) +v_{n}\left( r\right) \right) ^{\prime }}{%
\left( f_{1}+f_{2}\right) \left( u_{n}\left( r\right) +v_{n}\left( r\right)
,u_{n}\left( r\right) +v_{n}\left( r\right) \right) }\leq
r^{1-N}\int_{0}^{r}t^{N-1}\left( p_{1}\left( t\right) +p_{2}\left( t\right)
\right) dt,  \label{mat}
\end{equation}%
and%
\begin{equation*}
\int_{a_{1}+a_{2}}^{u_{n}\left( r\right) +v_{n}\left( r\right) }\frac{1}{%
f_{1}\left( t,t\right) +f_{2}\left( t,t\right) }dt\leq P_{1}\left( r\right)
+P_{2}\left( r\right) .
\end{equation*}%
We now have 
\begin{equation}
Z\left( u_{n}\left( r\right) +v_{n}\left( r\right) \right) \leq P_{1}\left(
r\right) +P_{2}\left( r\right) ,  \label{zc1}
\end{equation}%
which will play a basic role in the proof of our main results. The
inequalities (\ref{zc1}) can be rewritten as 
\begin{equation}
u_{n}\left( r\right) +v_{n}\left( r\right) \leq Z^{-1}\left( P_{1}\left(
r\right) +P_{2}\left( r\right) \right) .  \label{zc2}
\end{equation}%
This can be easily seen from the fact that $Z$ is a bijection with the
inverse function $Z$ strictly increasing on $\left[ 0,Z\left( \infty \right)
\right) $. Let $M_{1}\geq \max \left\{ 1,\frac{1}{a_{1}}\right\} $ and $%
M_{2}\geq \max \left\{ 1,\frac{1}{a_{2}}\right\} $. Then, going back to (\ref%
{gen1}) we have%
\begin{eqnarray}
\left[ r^{N-1}\left( u_{n}\left( r\right) \right) ^{\prime }\right] ^{\prime
} &=&r^{N-1}p_{1}\left( r\right) f_{1}\left( u_{n-1}\left( r\right)
,v_{n-1}\left( r\right) \right) \leq r^{N-1}p_{1}\left( r\right) f_{1}\left(
u_{n}\left( r\right) ,v_{n}\left( r\right) \right) ,  \notag \\
&\leq &r^{N-1}p_{1}\left( r\right) f_{1}\left( u_{n}\left( r\right)
,2u_{n}\left( r\right) +v_{n}\left( r\right) \right)  \notag \\
&\leq &r^{N-1}p_{1}\left( r\right) f_{1}\left( u_{n}\left( r\right)
,u_{n}\left( r\right) +Z^{-1}\left( P_{1}\left( r\right) +P_{2}\left(
r\right) \right) \right)  \notag \\
&\leq &r^{N-1}p_{1}\left( r\right) f_{1}\left( u_{n}\left( r\right)
,u_{n}\left( r\right) (1+\frac{1}{u_{n}\left( r\right) }Z^{-1}\left(
P_{1}\left( r\right) +P_{2}\left( r\right) \right) )\right)  \label{24} \\
&\leq &r^{N-1}p_{1}\left( r\right) f_{1}\left( u_{n}\left( r\right)
,u_{n}\left( r\right) (1+\frac{1}{a_{1}}Z^{-1}\left( P_{1}\left( r\right)
+P_{2}\left( r\right) \right) )\right)  \notag \\
&\leq &r^{N-1}p_{1}\left( r\right) f_{1}\left( u_{n}\left( r\right)
,u_{n}\left( r\right) M_{1}\left( 1+Z^{-1}\left( P_{1}\left( r\right)
+P_{2}\left( r\right) \right) \right) )\right)  \notag \\
&\leq &r^{N-1}p_{1}\left( r\right) \overline{c}_{1}f_{1}\left( u_{n}\left(
r\right) ,u_{n}\left( r\right) \right) \cdot \overline{f_{1}}\left(
M_{1}\left( 1+Z^{-1}\left( P_{1}\left( r\right) +P_{2}\left( r\right)
\right) \right) )\right) ,  \notag
\end{eqnarray}%
and in the same vein%
\begin{equation}
\left[ r^{N-1}\left( v_{n}\left( r\right) \right) ^{\prime }\right] ^{\prime
}\leq r^{N-1}p_{2}\left( r\right) \overline{c}_{2}f_{2}\left( v_{n}\left(
r\right) ,v_{n}\left( r\right) \right) \cdot \overline{f_{2}}\left(
M_{2}\left( 1+Z^{-1}\left( P_{1}\left( r\right) +P_{2}\left( r\right)
\right) \right) )\right) .  \label{25}
\end{equation}%
By (\ref{24}) and (\ref{25}), we have%
\begin{equation*}
\left\{ 
\begin{array}{l}
r^{N-1}\left( u_{n}\left( r\right) \right) ^{\prime \prime }\leq \left(
N-1\right) r^{N-2}\left( u_{n}\left( r\right) \right) ^{\prime
}+r^{N-1}\left( u_{n}\left( r\right) \right) ^{\prime \prime }=\left[
r^{N-1}\left( u_{n}\left( r\right) \right) ^{\prime }\right] ^{\prime } \\ 
\leq r^{N-1}p_{1}\left( r\right) \overline{c}_{1}f_{1}\left( u_{n}\left(
r\right) ,u_{n}\left( r\right) \right) \cdot \overline{f_{1}}\left(
M_{1}\left( 1+Z^{-1}\left( P_{1}\left( r\right) +P_{2}\left( r\right)
\right) \right) )\right) , \\ 
r^{N-1}\left( v_{n}\left( r\right) \right) ^{\prime \prime }\leq \left[
r^{N-1}\left( v_{n}\left( r\right) \right) ^{\prime }\right] ^{\prime } \\ 
\leq r^{N-1}p_{2}\left( r\right) \overline{c}_{2}f_{2}\left( v_{n}\left(
r\right) ,v_{n}\left( r\right) \right) \cdot \overline{f_{2}}\left(
M_{2}\left( 1+Z^{-1}\left( P_{1}\left( r\right) +P_{2}\left( r\right)
\right) \right) )\right) .%
\end{array}%
\right.
\end{equation*}%
Multiplying the first inequality by $\left( u_{n}\left( r\right) \right)
^{\prime }$ and the second by $\left( v_{n}\left( r\right) \right) ^{\prime
} $, we obtain%
\begin{equation}
\left\{ 
\begin{array}{l}
\left\{ \left[ \left( u_{n}\left( r\right) \right) ^{\prime }\right]
^{2}\right\} ^{\prime }\leq 2p_{1}\left( r\right) \overline{c}%
_{1}f_{1}\left( u_{n}\left( r\right) ,u_{n}\left( r\right) \right) \left(
u_{n}\left( r\right) \right) ^{\prime }\overline{f_{1}}\left( M_{1}\left(
1+Z^{-1}\left( P_{1}\left( r\right) +P_{2}\left( r\right) \right) \right)
)\right) , \\ 
\left\{ \left[ \left( v_{n}\left( r\right) \right) ^{\prime }\right]
^{2}\right\} ^{\prime }\leq 2p_{2}\left( r\right) \overline{c}%
_{2}f_{2}\left( v_{n}\left( r\right) ,v_{n}\left( r\right) \right) \left(
v_{n}\left( r\right) \right) ^{\prime }\overline{f_{2}}\left( M_{2}\left(
1+Z^{-1}\left( P_{1}\left( r\right) +P_{2}\left( r\right) \right) \right)
)\right) .%
\end{array}%
\right.  \label{27}
\end{equation}%
Integrating in (\ref{27}) from $0$ to $r$ we also have%
\begin{equation*}
\left\{ 
\begin{array}{l}
\left[ \left( u_{n}\left( r\right) \right) ^{\prime }\right] ^{2}\leq
2\int_{0}^{r}p_{1}\left( t\right) \overline{c}_{1}f_{1}\left( u_{n}\left(
t\right) ,u_{n}\left( t\right) \right) \left( u_{n}\left( t\right) \right)
^{\prime }\overline{f_{1}}\left( M_{1}\left( 1+Z^{-1}\left( P_{1}\left(
t\right) +P_{2}\left( t\right) \right) \right) )\right) dt, \\ 
\left[ \left( v_{n}\left( r\right) \right) ^{\prime }\right] ^{2}\leq
2\int_{0}^{r}p_{2}\left( t\right) \overline{c}_{2}f_{2}\left( v_{n}\left(
t\right) ,v_{n}\left( t\right) \right) \left( v_{n}\left( t\right) \right)
^{\prime }\overline{f_{2}}\left( M_{2}\left( 1+Z^{-1}\left( P_{1}\left(
t\right) +P_{2}\left( t\right) \right) \right) )\right) dt.%
\end{array}%
\right.
\end{equation*}%
Hence%
\begin{equation}
\left\{ 
\begin{array}{l}
\left[ \left( u_{n}\left( r\right) \right) ^{\prime }\right] ^{2}\leq 2%
\overline{f_{1}}\left( M_{1}\left( 1+Z^{-1}\left( P_{1}\left( r\right)
+P_{2}\left( r\right) \right) \right) )\right) \int_{0}^{r}p_{1}\left(
t\right) \overline{c}_{1}f_{1}\left( u_{n}\left( t\right) ,u_{n}\left(
t\right) \right) \left( u_{n}\left( t\right) \right) ^{\prime }dt, \\ 
\left[ \left( v_{n}\left( r\right) \right) ^{\prime }\right] ^{2}\leq 2%
\overline{f_{2}}\left( M_{2}\left( 1+Z^{-1}\left( P_{1}\left( r\right)
+P_{2}\left( r\right) \right) \right) )\right) \int_{0}^{r}p_{2}\left(
t\right) \overline{c}_{2}f_{2}\left( v_{n}\left( t\right) ,v_{n}\left(
t\right) \right) \left( v_{n}\left( t\right) \right) ^{\prime }dt.%
\end{array}%
\right.  \label{ineqr}
\end{equation}%
Set now 
\begin{eqnarray}
\phi _{1}\left( r\right) &=&\max \left\{ p_{1}\left( t\right) \left\vert
0\leq t\leq r\right. \right\} \text{,}  \label{mom2} \\
\phi _{2}\left( r\right) &=&\max \left\{ p_{2}\left( t\right) \left\vert
0\leq t\leq r\right. \right\} \text{.}  \notag
\end{eqnarray}%
Thanks to the definition of $\phi _{1}\left( r\right) $ and $\phi _{2}\left(
r\right) $ we get from the inequalities (\ref{ineqr}) that%
\begin{equation}
\left\{ 
\begin{array}{l}
\left[ \left( u_{n}\left( r\right) \right) ^{\prime }\right] ^{2}\leq 2%
\overline{c}_{1}\overline{f_{1}}\left( M_{1}\left( 1+Z^{-1}\left(
P_{1}\left( r\right) +P_{2}\left( r\right) \right) \right) )\right) \phi
_{1}\left( r\right) \int_{0}^{r}f_{1}\left( u_{n}\left( t\right)
,u_{n}\left( t\right) \right) \left( u_{n}\left( t\right) \right) ^{\prime
}dt, \\ 
\left[ \left( v_{n}\left( r\right) \right) ^{\prime }\right] ^{2}\leq 2%
\overline{c}_{2}\overline{f_{2}}\left( M_{2}\left( 1+Z^{-1}\left(
P_{1}\left( r\right) +P_{2}\left( r\right) \right) \right) )\right) \phi
_{2}\left( r\right) \int_{0}^{r}f_{2}\left( v_{n}\left( t\right)
,v_{n}\left( t\right) \right) \left( v_{n}\left( t\right) \right) ^{\prime
}dt.%
\end{array}%
\right.  \label{210}
\end{equation}%
As a consequence of (\ref{210}), we also have%
\begin{equation}
\left\{ 
\begin{array}{l}
\left( u_{n}\left( r\right) \right) ^{\prime }\leq \sqrt{2\overline{c}_{1}%
\overline{f_{1}}\left( M_{1}\left( 1+Z^{-1}\left( P_{1}\left( r\right)
+P_{2}\left( r\right) \right) \right) )\right) \phi _{1}\left( r\right) }%
\sqrt{\int_{a_{1}}^{u_{n}\left( r\right) }f_{1}\left( t,t\right) dt}, \\ 
\left( v_{n}\left( r\right) \right) ^{\prime }\leq \sqrt{2\overline{c}_{2}%
\overline{f_{2}}\left( M_{2}\left( 1+Z^{-1}\left( P_{1}\left( r\right)
+P_{2}\left( r\right) \right) \right) )\right) \phi _{2}\left( r\right) }%
\sqrt{\int_{a_{2}}^{v_{n}\left( r\right) }f_{2}\left( t,t\right) dt},%
\end{array}%
\right.  \label{211}
\end{equation}%
and, thus%
\begin{equation}
\left\{ 
\begin{array}{l}
\frac{\left( u_{n}\left( r\right) \right) ^{\prime }}{\sqrt{%
\int_{a_{1}}^{u_{n}\left( r\right) }f_{1}\left( t,t\right) dt}}\leq \sqrt{2%
\overline{c}_{1}\overline{f_{1}}\left( M_{1}\left( 1+Z^{-1}\left(
P_{1}\left( r\right) +P_{2}\left( r\right) \right) \right) )\right) \phi
_{1}\left( r\right) }, \\ 
\frac{\left( v_{n}\left( r\right) \right) ^{\prime }}{\sqrt{%
\int_{a_{2}}^{v_{n}\left( r\right) }f_{2}\left( t,t\right) dt}}\leq \sqrt{2%
\overline{c}_{2}\overline{f_{2}}\left( M_{2}\left( 1+Z^{-1}\left(
P_{1}\left( r\right) +P_{2}\left( r\right) \right) \right) )\right) \phi
_{2}\left( r\right) }.%
\end{array}%
\right.  \label{212}
\end{equation}%
Integrating (\ref{212}) leads to%
\begin{equation*}
\left\{ 
\begin{array}{l}
\int_{a_{1}}^{u_{n}\left( r\right) }\left( \int_{0}^{s}f_{1}\left(
t,t\right) dt\right) ^{-1/2}ds\leq \sqrt{2\overline{c}_{1}\overline{f_{1}}%
\left( M_{1}\left( 1+Z^{-1}\left( P_{1}\left( r\right) +P_{2}\left( r\right)
\right) \right) )\right) }\int_{0}^{r}\sqrt{\phi _{1}\left( s\right) }ds, \\ 
\int_{a_{2}}^{v_{n}\left( r\right) }\left( \int_{a_{2}}^{s}f_{2}\left(
t,t\right) dt\right) ^{-1/2}ds\leq \sqrt{2\overline{c}_{2}\overline{f_{2}}%
\left( M_{2}\left( 1+Z^{-1}\left( P_{1}\left( r\right) +P_{2}\left( r\right)
\right) \right) )\right) }\int_{0}^{r}\sqrt{\phi _{2}\left( s\right) }ds.%
\end{array}%
\right.
\end{equation*}%
which is the same as%
\begin{equation}
\mathcal{KO}_{f_{1}}\left( u_{n}\left( r\right) \right) \leq \sqrt{2%
\overline{c}_{1}}\cdot \overline{P}_{1}\left( r\right) \text{ and }\mathcal{%
KO}_{f_{2}}\left( v_{n}\left( r\right) \right) \leq \sqrt{2\overline{c}_{2}}%
\cdot \overline{P}_{2}\left( r\right) .  \label{ints}
\end{equation}%
Now, we can easy see that $\mathcal{KO}_{f_{1}}$ is a bijection with the
inverse function $\mathcal{KO}_{f_{1}}^{-1}$ strictly increasing on $\left[
0,\mathcal{KO}_{f_{1}}\left( \infty \right) \right) $. Hence%
\begin{equation}
u_{n}\left( r\right) \leq \mathcal{KO}_{f_{1}}^{-1}\left( \sqrt{2\overline{c}%
_{1}}\cdot \overline{P}_{1}\left( r\right) \right)  \label{int}
\end{equation}%
The second inequality in (\ref{ints}), leads to%
\begin{equation}
v_{n}\left( r\right) \leq \mathcal{KO}_{f_{2}}^{-1}\left( \sqrt{2\overline{c}%
_{2}}\cdot \overline{P}_{2}\left( r\right) \right) .  \label{int2}
\end{equation}%
In summary, we have found upper bounds for $\left\{ u_{n}\right\} _{n\geq 0}$
and $\left\{ v_{n}\right\} _{n\geq 0}$ which are dependent of $r$. Now let
us complete the proof of Theorems \ref{th1}-\ref{th5}. We prove that the
sequences $\left\{ u_{n}\right\} _{n\geq 0}$ and $\left\{ v_{n}\right\}
_{n\geq 0}$ are bounded and equicontinuous on $\left[ 0,c_{0}\right] $ for
arbitrary $c_{0}>0$. Indeed, since 
\begin{equation*}
\left( u_{n}\left( r\right) \right) ^{^{\prime }}\geq 0\text{ and }\left(
v_{n}\left( r\right) \right) ^{^{\prime }}\geq 0\text{ for all }r\geq 0,
\end{equation*}%
it follows that 
\begin{equation*}
u_{n}\left( r\right) \leq u_{n}\left( c_{0}\right) \leq C_{1}\text{ and }%
v_{n}\left( r\right) \leq v_{n}\left( c_{0}\right) \leq C_{2}\text{ on }%
\left[ 0,c_{0}\right] .
\end{equation*}%
Here $C_{1}=\mathcal{KO}_{f_{1}}^{-1}\left( \sqrt{2\overline{c}_{1}}\cdot 
\overline{P}_{1}\left( c_{0}\right) \right) $ and $C_{2}=\mathcal{KO}%
_{f_{2}}^{-1}\left( \sqrt{2\overline{c}_{2}}\cdot \overline{P}_{2}\left(
c_{0}\right) \right) $ are positive constants. Recall that $\left\{
u_{n}\right\} _{n\geq 0}$ and $\left\{ v_{n}\right\} _{n\geq 0}$ are bounded
on $\left[ 0,c_{0}\right] $ for arbitrary $c_{0}>0$. Using this fact, we
show that the same is true of $\left( u_{n}\left( r\right) \right) ^{\prime
} $ and $\left( v_{n}\left( r\right) \right) ^{\prime }$. Indeed, for any $%
r\geq 0$,%
\begin{eqnarray*}
\left( u_{n}\left( r\right) \right) ^{\prime }
&=&r^{1-N}\int_{0}^{r}t^{N-1}p_{1}\left( t\right) f_{1}\left( u_{n-1}\left(
t\right) ,v_{n-1}\left( t\right) \right) dt \\
&\leq &r^{1-N}r^{N-1}\int_{0}^{r}p_{1}\left( t\right) f_{1}\left(
u_{n}\left( t\right) ,v_{n}\left( t\right) \right) dt \\
&\leq &\left\Vert p_{1}\right\Vert _{\infty }\int_{0}^{r}f_{1}\left(
u_{n}\left( t\right) ,v_{n}\left( t\right) \right) dt\leq \left\Vert
p_{1}\right\Vert _{\infty }f_{1}\left( C_{1},C_{2}\right) \int_{0}^{r}dt \\
&\leq &\left\Vert p_{1}\right\Vert _{\infty }f_{1}\left( C_{1},C_{2}\right)
c_{0}\text{ on }\left[ 0,c_{0}\right] .
\end{eqnarray*}%
Similar arguments show that%
\begin{eqnarray*}
\left( v_{n}\left( r\right) \right) ^{\prime }
&=&r^{1-N}\int_{0}^{r}t^{N-1}p_{2}\left( t\right) f_{2}\left( u_{n-1}\left(
t\right) ,v_{n-1}\left( t\right) \right) dt \\
&\leq &r^{1-N}f_{2}\left( C_{1},C_{2}\right) r^{N-1}\left\Vert
p_{2}\right\Vert _{\infty }\int_{0}^{r}dt \\
&\leq &\left\Vert p_{2}\right\Vert _{\infty }f_{2}\left( C_{1},C_{2}\right)
c_{0}\text{ on }\left[ 0,c_{0}\right] .
\end{eqnarray*}%
It remains, to prove that $\left\{ u_{n}\right\} _{n\geq 0}$ and $\left\{
v_{n}\right\} _{n\geq 0}$ are equicontinuous on $\left[ 0,c_{0}\right] $ for
arbitrary $c_{0}>0$. Let $\varepsilon _{1}$, $\varepsilon _{2}>0$. To verify
equicontinuous on $\left[ 0,c_{0}\right] $, observe that%
\begin{eqnarray*}
\left\vert u_{n}\left( x\right) -u_{n}\left( y\right) \right\vert
&=&\left\vert \left( u_{n}\left( \xi _{1}\right) \right) ^{\prime
}\right\vert \left\vert x-y\right\vert \leq \left\Vert p_{1}\right\Vert
_{\infty }f_{1}\left( C_{1},C_{2}\right) c_{0}\left\vert x-y\right\vert , \\
\left\vert v_{n}\left( x\right) -v_{n}\left( y\right) \right\vert
&=&\left\vert \left( v_{n}\left( \xi _{2}\right) \right) ^{\prime
}\right\vert \left\vert x-y\right\vert \leq \left\Vert p_{2}\right\Vert
_{\infty }f_{2}\left( C_{1},C_{2}\right) c_{0}\left\vert x-y\right\vert ,
\end{eqnarray*}%
for all $n\in \mathbb{N}$ and all $x,y\in \left[ 0,c_{0}\right] $ and for $%
\xi _{1}$, $\xi _{2}$ the constants from the mean value theorem. So it
suffices to take%
\begin{equation*}
\delta _{1}=\frac{\varepsilon _{1}}{\left\Vert p_{1}\right\Vert _{\infty
}f_{1}\left( C_{1},C_{2}\right) c_{0}}\text{ and }\delta _{2}=\frac{%
\varepsilon _{2}}{f_{2}\left( C_{1},C_{2}\right) \left\Vert p_{2}\right\Vert
_{\infty }c_{0}},
\end{equation*}%
to see that $\left\{ u_{n}\right\} _{n\geq 0}$ and $\left\{ v_{n}\right\}
_{n\geq 0}$ are equicontinuous on $\left[ 0,c_{0}\right] $. Thus, it follows
from the Arzela--Ascoli theorem that there exists a function $u\in C\left( %
\left[ 0,c_{0}\right] \right) $ and a subsequence $N_{1}$ of $\mathbb{N}%
^{\ast }$ with $u_{n}\left( r\right) $ converging uniformly on $u$ to $\left[
0,c_{0}\right] $ as $n\rightarrow \infty $ through $N_{1}$. By the same
token there exists a function $v\in C\left( \left[ 0,c_{0}\right] \right) $
and a subsequence $N_{2}$ of $\mathbb{N}^{\ast }$ with $v_{n}\left( r\right) 
$ converging uniformly to $v$ on $\left[ 0,c_{0}\right] $ as $n\rightarrow
\infty $ through $N_{2}$. Thus $\left\{ \left( u_{n}\left( r\right)
,v_{n}\left( r\right) \right) \right\} _{n\in N_{2}}$ converges uniformly on 
$\left[ 0,c_{0}\right] $ to $\left( u,v\right) \in C\left( \left[ 0,c_{0}%
\right] \right) \times C\left( \left[ 0,c_{0}\right] \right) $ through $%
N_{2} $ (see L\"{u}-O'Regan-Agarwal \cite{LU}). The solution constructed in
this way will be radially symmetric solution of system (\ref{11}). We remark
that, since $\left\{ \left( u_{n}\left( r\right) ,v_{n}\left( r\right)
\right) \right\} _{n}$ is non-decreasing on $[0,\infty )$, then $\left\{
\left( u_{n}\left( r\right) ,v_{n}\left( r\right) \right) \right\} _{n}$
itself converges uniformly to $\left( u,v\right) $ on $[0,c_{0}]$. Moreover,
the radial solutions of (\ref{11}) with $u\left( 0\right) =a_{1},$ $v\left(
0\right) =a_{2}$ satisfy:%
\begin{eqnarray}
u\left( r\right) &=&a_{1}+\int_{0}^{r}y^{1-N}\int_{0}^{y}t^{N-1}p_{1}\left(
t\right) f_{1}\left( u\left( t\right) ,v\left( t\right) \right) dtdy,\text{ }%
r\geq 0,  \label{eq1} \\
v\left( r\right) &=&a_{2}+\int_{0}^{r}y^{1-N}\int_{0}^{y}t^{N-1}p_{2}\left(
t\right) f_{2}\left( u\left( t\right) ,v\left( t\right) \right) dtdy,\text{ }%
r\geq 0.  \label{eq2}
\end{eqnarray}%
It remains to show that $u\in C^{2}\left[ 0,\infty \right) $. It is not
difficult to see that $u\in C^{2}\left( 0,\infty \right) \cap C\left[
0,\infty \right) $. We prove that $u^{\prime }\left( r\right) $ and $%
u^{\prime \prime }\left( r\right) $ are continuous at $r=0$. To ensure that $%
u^{\prime }\left( r\right) $ is continuous at $r=0$ we evaluate%
\begin{equation*}
\lim_{r\rightarrow 0}u^{\prime }\left( r\right) =\lim_{r\rightarrow 0}\frac{%
\int_{0}^{r}s^{N-1}p_{1}\left( s\right) f_{1}\left( u\left( s\right)
,v\left( s\right) \right) ds}{r^{N-1}}.
\end{equation*}%
By applying L'Hopital's rule, we have%
\begin{equation*}
\lim_{r\rightarrow 0}u^{\prime }\left( r\right) =\lim_{r\rightarrow 0}\frac{%
r^{N-1}p_{1}\left( r\right) f_{1}\left( u\left( r\right) ,v\left( r\right)
\right) }{\left( N-1\right) r^{N-2}}=\lim_{r\rightarrow 0}\frac{rp_{1}\left(
r\right) f_{1}\left( u\left( r\right) ,v\left( r\right) \right) }{\left(
N-1\right) }=0.
\end{equation*}%
On the other hand by a successive application of the L'Hopital's rule, we
have%
\begin{eqnarray*}
u^{\prime }\left( 0\right) &=&\lim_{r\rightarrow 0}\frac{u\left( r\right)
-u\left( 0\right) }{r} \\
&=&\lim_{r\rightarrow 0}\frac{\int_{0}^{r}y^{1-N}\int_{0}^{y}t^{N-1}p_{1}%
\left( t\right) f_{1}\left( u\left( t\right) ,v\left( t\right) \right) dtdy}{%
r} \\
&=&\lim_{r\rightarrow 0}\frac{\int_{0}^{r}t^{N-1}p_{1}\left( t\right)
f_{1}\left( u\left( t\right) ,v\left( t\right) \right) dt}{r^{N-1}} \\
&=&\lim_{r\rightarrow 0}\frac{r^{N-1}p_{1}\left( r\right) f_{1}\left(
u\left( r\right) ,v\left( r\right) \right) }{N-1}=0.
\end{eqnarray*}%
Consequently, 
\begin{equation*}
\lim_{r\rightarrow 0}u^{\prime }\left( r\right) =u^{\prime }\left( 0\right)
=0
\end{equation*}%
and the proof of $u^{\prime }\left( r\right) $ is continuous at $r=0$ is
concluded. We now prove that $u^{\prime \prime }\left( r\right) $ is
continuous at $r=0$. To this aim, remark that 
\begin{equation*}
u^{\prime \prime }\left( r\right) =\left( 1-N\right)
r^{-N}\int_{0}^{r}s^{N-1}p_{1}\left( s\right) f_{1}\left( u\left( s\right)
,v\left( s\right) \right) ds+p_{1}\left( r\right) f_{1}\left( u\left(
r\right) ,v\left( r\right) \right) \text{, }r>0.
\end{equation*}%
We deduce%
\begin{equation*}
\lim_{r\rightarrow 0}u^{\prime \prime }\left( r\right) =\lim_{r\rightarrow 0}%
\frac{\left( 1-N\right) \int_{0}^{r}s^{N-1}p_{1}\left( s\right) f_{1}\left(
u\left( s\right) ,v\left( s\right) \right) ds}{r^{N}}+\lim_{r\rightarrow
0}p_{1}\left( r\right) f_{1}\left( u\left( r\right) ,v\left( r\right)
\right) .
\end{equation*}%
In the sequel, we use the following notation%
\begin{equation*}
L_{1}\left( r\right) =\frac{\left( 1-N\right) \int_{0}^{r}s^{N-1}p_{1}\left(
s\right) f_{1}\left( u\left( s\right) ,v\left( s\right) \right) ds}{r^{N}}%
\text{ and }L_{2}\left( r\right) =p_{1}\left( r\right) f_{1}\left( u\left(
r\right) ,v\left( r\right) \right) .
\end{equation*}%
According to our notation, we remark that 
\begin{eqnarray*}
\lim_{r\rightarrow 0}L_{1}\left( r\right) &=&\lim_{r\rightarrow 0}\frac{%
\left( 1-N\right) r^{N-1}p_{1}\left( r\right) f_{1}\left( u\left( r\right)
,v\left( r\right) \right) }{Nr^{N-1}} \\
&=&\lim_{r\rightarrow 0}\left( \frac{p_{1}\left( r\right) f_{1}\left(
u\left( r\right) ,v\left( r\right) \right) }{N}-p_{1}\left( r\right)
f_{1}\left( u\left( r\right) ,v\left( r\right) \right) \right) \\
&=&\frac{p_{1}\left( 0\right) f_{1}\left( u\left( 0\right) ,v\left( 0\right)
\right) }{N}-p_{1}\left( 0\right) f_{1}\left( u\left( 0\right) ,v\left(
0\right) \right) \\
\lim_{r\rightarrow 0}L_{2}\left( r\right) &=&p_{1}\left( 0\right)
f_{1}\left( u\left( 0\right) ,v\left( 0\right) \right) .
\end{eqnarray*}%
Clearly%
\begin{equation}
\lim_{r\rightarrow 0}u^{\prime \prime }\left( r\right) =\lim_{r\rightarrow
0}L_{1}\left( r\right) +\lim_{r\rightarrow 0}L_{2}\left( r\right) =\frac{%
p_{1}\left( 0\right) f_{1}\left( u\left( 0\right) ,v\left( 0\right) \right) 
}{N}.
\end{equation}%
On the other hand%
\begin{equation*}
u^{\prime \prime }\left( 0\right) =\lim_{r\rightarrow 0}\frac{u^{\prime
}\left( r\right) -u^{\prime }\left( 0\right) }{r}=\lim_{r\rightarrow 0}\frac{%
r^{1-N}\int_{0}^{r}s^{N-1}p_{1}\left( s\right) f_{1}\left( u\left( s\right)
,v\left( s\right) \right) ds}{r},
\end{equation*}%
where we have used $u^{\prime }\left( 0\right) =0$. Now, by the use of
L'Hopital's rule%
\begin{equation*}
u^{\prime \prime }\left( r\right) =\lim_{r\rightarrow 0}\left( 1-N\right)
r^{-N}\int_{0}^{r}s^{N-1}p_{1}\left( s\right) f_{1}\left( u\left( s\right)
,v\left( s\right) \right) ds+\lim_{r\rightarrow 0}p_{1}\left( r\right)
f_{1}\left( u\left( r\right) ,v\left( r\right) \right) .
\end{equation*}%
Rearranging, we get%
\begin{equation*}
u^{\prime \prime }\left( 0\right) =\lim_{r\rightarrow 0}\left( L_{1}\left(
r\right) +L_{2}\left( r\right) \right) .
\end{equation*}%
It remains now to see that 
\begin{equation}
u^{\prime \prime }\left( 0\right) =\frac{p_{1}\left( 0\right) f_{1}\left(
u\left( 0\right) ,v\left( 0\right) \right) }{N}.
\end{equation}%
Then $u^{\prime \prime }$ is continuous at $r=0$. Hence the claim follows
from (\ref{210}) and (\ref{211}). This means that $u\in C^{2}\left[ 0,\infty
\right) $. In the same vein $v\in C^{2}\left[ 0,\infty \right) $. The
existence and regularity of solutions for the Theorems \ref{th1}-\ref{th5}
was proved. Next, choose $R>0$ so that $r^{2N-2}p_{1}\left( r\right) $ and $%
r^{2N-2}p_{2}\left( r\right) $ are non-decreasing for $r\geq R$. In order,
to prove Theorems \ref{th1}-\ref{th5} we intend to establish some
inequalities. Using the same arguments as in (\ref{24}) and (\ref{25}) we
can see that%
\begin{equation}
\left\{ 
\begin{array}{l}
\left[ r^{N-1}\left( u\left( r\right) \right) ^{\prime }\right] ^{\prime
}\leq r^{N-1}p_{1}\left( r\right) \overline{c}_{1}f_{1}\left( u\left(
r\right) ,u\left( r\right) \right) \overline{f_{1}}\left( M_{1}\left(
1+Z^{-1}\left( P_{1}\left( r\right) +P_{2}\left( r\right) \right) \right)
)\right) , \\ 
\left[ r^{N-1}\left( v\left( r\right) \right) ^{\prime }\right] ^{\prime
}\leq r^{N-1}p_{2}\left( r\right) \overline{c}_{2}f_{2}\left( v\left(
r\right) ,v\left( r\right) \right) \overline{f_{2}}\left( M_{2}\left(
1+Z^{-1}\left( P_{1}\left( r\right) +P_{2}\left( r\right) \right) \right)
)\right) .%
\end{array}%
\right.  \label{217}
\end{equation}%
Multiplying the first equation in (\ref{217}) by $r^{N-1}\left( u\left(
r\right) \right) ^{\prime }$ and the second by $r^{N-1}\left( v\left(
r\right) \right) ^{\prime }$ and integrating gives%
\begin{equation*}
\left\{ 
\begin{array}{l}
\left[ r^{N-1}\left( u\left( r\right) \right) ^{\prime }\right] ^{2}\leq %
\left[ R^{N-1}\left( u\left( R\right) \right) ^{\prime }\right] ^{2} \\ 
+2\int_{R}^{r}t^{2\left( N-1\right) }p_{1}\left( t\right) \overline{c}_{1}%
\overline{f_{1}}\left( M_{1}\left( 1+Z^{-1}\left( P_{1}\left( t\right)
+P_{2}\left( t\right) \right) \right) )\right) \frac{d}{dz}%
\int_{a_{1}}^{u\left( z\right) }f_{1}\left( s,s\right) dsdt, \\ 
\left[ r^{N-1}\left( v\left( r\right) \right) ^{\prime }\right] ^{2}\leq %
\left[ R^{N-1}\left( v\left( R\right) \right) ^{\prime }\right] ^{2} \\ 
+2\int_{R}^{r}t^{2\left( N-1\right) }p_{2}\left( t\right) \overline{c}_{2}%
\overline{f_{2}}\left( M_{2}\left( 1+Z^{-1}\left( P_{1}\left( t\right)
+P_{2}\left( t\right) \right) \right) )\right) \frac{d}{dz}%
\int_{a_{2}}^{v\left( z\right) }f_{2}\left( s,s\right) dsdt,%
\end{array}%
\right.
\end{equation*}%
for $r\geq R$. We get from the monotonicity of $z^{2N-2}p_{1}\left( z\right) 
$ and $z^{2N-2}p_{2}\left( z\right) $ for $r\geq z\geq R$ that%
\begin{equation*}
\left\{ 
\begin{array}{l}
\left[ r^{N-1}\left( u\left( r\right) \right) ^{\prime }\right] ^{2}\leq
C_{1}+2r^{N-1}p_{1}\left( r\right) \overline{c}_{1}\overline{f_{1}}\left(
M_{1}\left( 1+Z^{-1}\left( P_{1}\left( r\right) +P_{2}\left( r\right)
\right) \right) )\right) \int_{0}^{u\left( r\right) }f_{1}\left( z,z\right)
dz, \\ 
\left[ r^{N-1}\left( v\left( r\right) \right) ^{\prime }\right] ^{2}\leq
C_{2}+2r^{N-1}p_{2}\left( r\right) \overline{c}_{2}\overline{f_{2}}\left(
M_{2}\left( 1+Z^{-1}\left( P_{1}\left( r\right) +P_{2}\left( r\right)
\right) \right) )\right) \int_{0}^{v\left( r\right) }f_{2}\left( z,z\right)
dz,%
\end{array}%
\right.
\end{equation*}%
where $C_{1}=\left[ R^{N-1}\left( u\left( R\right) \right) ^{\prime }\right]
^{2}$ and $C_{2}=\left[ R^{N-1}\left( v\left( R\right) \right) ^{\prime }%
\right] ^{2}$. This implies that%
\begin{equation}
\left\{ 
\begin{array}{l}
\frac{\left( u\left( r\right) \right) ^{\prime }}{\sqrt{\int_{0}^{u\left(
r\right) }f_{1}\left( z,z\right) dz}}\leq \frac{r^{1-N}\sqrt{C_{1}}}{\sqrt{%
\int_{0}^{u\left( r\right) }f_{1}\left( z,z\right) dz}}+\sqrt{2p_{1}\left(
r\right) \overline{c}_{1}\overline{f_{1}}\left( M_{1}\left( 1+Z^{-1}\left(
P_{1}\left( r\right) +P_{2}\left( r\right) \right) \right) )\right) }, \\ 
\frac{\left( v\left( r\right) \right) ^{\prime }}{\sqrt{\int_{0}^{v\left(
r\right) }f_{2}\left( z,z\right) dz}}\leq \frac{r^{1-N}\sqrt{C_{2}}}{\sqrt{%
\int_{0}^{v\left( r\right) }f_{2}\left( z,z\right) dz}}+\sqrt{2p_{2}\left(
r\right) \overline{c}_{2}\overline{f_{2}}\left( M_{2}\left( 1+Z^{-1}\left(
P_{1}\left( r\right) +P_{2}\left( r\right) \right) \right) )\right) }.%
\end{array}%
\right.  \label{218}
\end{equation}%
In particular, integrating (\ref{218}) from $R$ tor $r$ and using the fact
that%
\begin{eqnarray*}
&&\sqrt{2r^{\frac{1+\varepsilon }{2}}p_{1}\left( r\right) \overline{c}_{1}%
\overline{f_{1}}\left( M_{1}\left( 1+Z^{-1}\left( P_{1}\left( r\right)
+P_{2}\left( r\right) \right) \right) )\right) r^{-\frac{1+\varepsilon }{2}}}
\\
&\leq &r^{1+\varepsilon }p_{1}\left( r\right) \overline{c}_{1}\overline{f_{1}%
}\left( M_{1}\left( 1+Z^{-1}\left( P_{1}\left( r\right) +P_{2}\left(
r\right) \right) \right) )\right) +r^{-\left( 1+\varepsilon \right) }
\end{eqnarray*}%
and%
\begin{eqnarray*}
&&\sqrt{2r^{\frac{1+\varepsilon }{2}}p_{2}\left( r\right) \overline{c}_{2}%
\overline{f_{2}}\left( M_{2}\left( 1+Z^{-1}\left( P_{1}\left( r\right)
+P_{2}\left( r\right) \right) \right) )\right) r^{-\frac{1+\varepsilon }{2}}}
\\
&\leq &r^{1+\varepsilon }p_{2}\left( r\right) \overline{c}_{2}\overline{f_{2}%
}\left( M_{2}\left( 1+Z^{-1}\left( P_{1}\left( r\right) +P_{2}\left(
r\right) \right) \right) )\right) +r^{-\left( 1+\varepsilon \right) }
\end{eqnarray*}%
for every $\varepsilon >0$, lead to%
\begin{eqnarray}
\int_{u\left( R\right) }^{u\left( r\right) }\left( \int_{0}^{t}f_{1}\left(
z,z\right) dz\right) ^{-1/2}dt &=&\mathcal{KO}_{f_{1}}\left( u\left(
r\right) \right) -\mathcal{KO}_{f_{1}}\left( u\left( R\right) \right)  \notag
\\
&\leq &\sqrt{C_{1}}\int_{R}^{r}\frac{t^{1-N}}{\sqrt{\int_{0}^{u\left(
t\right) }f_{1}\left( z,z\right) dz}}dt+\overline{c}_{1}\overline{P}%
_{1\varepsilon }+\int_{R}^{r}t^{-\left( 1+\varepsilon \right) }dt
\label{219} \\
&\leq &\sqrt{C_{1}}\frac{\int_{R}^{r}t^{1-N}dt}{\sqrt{\int_{0}^{u\left(
R\right) }f_{1}\left( z,z\right) dz}}+\overline{c}_{1}\overline{P}%
_{1\varepsilon }+\frac{1}{\varepsilon R^{\varepsilon }}.  \notag
\end{eqnarray}%
We next turn to estimating $v\left( r\right) $. A similar calculation yields%
\begin{equation}
\int_{v\left( R\right) }^{v\left( r\right) }\left( \int_{0}^{t}f_{2}\left(
z,z\right) dz\right) ^{-1/2}dt\leq \sqrt{C_{2}}\frac{\int_{R}^{r}t^{1-N}dt}{%
\sqrt{\int_{0}^{v\left( R\right) }f_{2}\left( z,z\right) dz}}+\overline{c}%
_{2}\overline{P}_{2\varepsilon }+\frac{1}{\varepsilon R^{\varepsilon }}.
\label{220}
\end{equation}%
The inequalities (\ref{219}) and (\ref{220}) are needed in proving the
\textquotedblright boundedness\textquotedblright\ of the functions $u$ and $%
v $. Indeed, they can be written as%
\begin{equation}
\left\{ 
\begin{array}{c}
u\left( r\right) \leq \mathcal{KO}_{f_{1}}^{-1}\left( \sqrt{C_{1}}%
\int_{R}^{r}t^{1-N}dt\left( \int_{0}^{u\left( R\right) }f_{1}\left(
z,z\right) dz\right) ^{-1/2}+\overline{c}_{1}\overline{P}_{1\varepsilon }+%
\frac{1}{\varepsilon R^{\varepsilon }}\right) , \\ 
v\left( r\right) \leq \mathcal{KO}_{f_{2}}^{-1}\left( \sqrt{C_{2}}%
\int_{R}^{r}t^{1-N}dt\left( \int_{0}^{v\left( R\right) }f_{2}\left(
z,z\right) dz\right) ^{-1/2}+\overline{c}_{2}\overline{P}_{2\varepsilon }+%
\frac{1}{\varepsilon R^{\varepsilon }}\right) .%
\end{array}%
\right.  \label{exp}
\end{equation}%
Next we prove that all Theorems hold.

\subsection{\textbf{Proof of Theorem \protect\ref{th1} completed:}}

In the case $\underline{P}\left( \infty \right) =\underline{Q}\left( \infty
\right) =\infty $, we observe that 
\begin{eqnarray}
u\left( r\right) &=&a_{1}+\int_{0}^{r}t^{1-N}\int_{0}^{t}s^{N-1}p_{1}\left(
s\right) f_{1}\left( u\left( s\right) ,v\left( s\right) \right) dsdt  \notag
\\
&\geq &a_{1}+\int_{0}^{r}y^{1-N}\int_{0}^{y}t^{N-1}p_{1}\left( t\right)
f_{1}\left( a_{1},a_{2}+f_{2}\left( a_{1},a_{2}\right) P_{2}\left( t\right)
\right) dtdy  \label{ints1} \\
&=&a_{1}+\underline{P}\left( r\right) .  \notag
\end{eqnarray}%
Proceeding as in the above, we also have 
\begin{eqnarray}
v\left( r\right) &=&a_{2}+\int_{0}^{r}y^{1-N}\int_{0}^{y}t^{N-1}p_{2}\left(
t\right) f_{2}\left( u\left( t\right) ,v\left( t\right) \right) dtdy  \notag
\\
&\geq &a_{2}+\int_{0}^{r}y^{1-N}\int_{0}^{y}t^{N-1}p_{2}\left( t\right)
f_{2}\left( a_{1},a_{1}+f_{1}\left( a_{1},a_{2}\right) P_{1}\left( t\right)
\right) dtdy  \label{ints2} \\
&=&a_{2}+\underline{Q}\left( r\right) .  \notag
\end{eqnarray}%
By taking limits in (\ref{ints1}) and (\ref{ints2}), we get entire large
solutions 
\begin{equation*}
\lim_{r\rightarrow \infty }u\left( r\right) =\infty \text{ and }%
\lim_{r\rightarrow \infty }v\left( r\right) =\infty .
\end{equation*}%
This completes the proof. We next consider:

\subsection{\textbf{Proof of Theorem \protect\ref{th2} completed:}}

When $\overline{P}_{1\varepsilon }\left( \infty \right) <\infty $ and $%
\overline{P}_{2\varepsilon }\left( \infty \right) <\infty $ we find from (%
\ref{exp}) that%
\begin{equation*}
u\left( r\right) <\infty \text{ and }v\left( r\right) <\infty \text{ for all 
}r\geq 0.
\end{equation*}%
In other words, we get that $\left( u,v\right) $ is bounded. \ Next, we give
the proof for the case where the solution is entire large. As in (\ref{219})
and (\ref{220}) we get 
\begin{equation}
\int_{u\left( R\right) }^{u\left( r\right) }\left( \int_{0}^{t}f_{1}\left(
z,z\right) dz\right) ^{-1/2}dt\leq C_{R}^{1}+\overline{c}_{1}\overline{P}%
_{1\varepsilon },  \label{f1}
\end{equation}%
and%
\begin{equation}
\int_{v\left( R\right) }^{v\left( r\right) }\left( \int_{0}^{t}f_{2}\left(
z,z\right) dz\right) ^{-1/2}dt\leq C_{R}^{2}+\overline{c}_{2}\overline{P}%
_{2\varepsilon },  \label{f2}
\end{equation}%
where 
\begin{equation*}
C_{R}^{1}=\sqrt{C_{1}}\frac{\int_{R}^{r}t^{1-N}dt}{\sqrt{\int_{0}^{u\left(
R\right) }f_{1}\left( z,z\right) dz}}+\frac{1}{\varepsilon R^{\varepsilon }}%
\text{ and }C_{R}^{2}=\sqrt{C_{2}}\frac{\int_{R}^{r}t^{1-N}dt}{\left(
\int_{0}^{v\left( R\right) }f_{2}\left( z,z\right) dz\right) ^{-1/2}}+\frac{1%
}{\varepsilon R^{\varepsilon }}.
\end{equation*}%
Letting $r\rightarrow \infty $ in (\ref{f1}) and (\ref{f2}), we find that $%
p_{1}$ and $p_{2}$ satisfy $\overline{P}_{1\varepsilon }=\infty $ and $%
\overline{P}_{2\varepsilon }=\infty $.

\subsection{\textbf{Proof of Theorem \protect\ref{th3} completed:}}

\textbf{1):} In the spirit of Theorem \ref{th1} and Theorem \ref{th2}, we
have 
\begin{equation*}
u\left( r\right) \leq \mathcal{KO}_{f_{1}}^{-1}\left( \sqrt{C_{1}}%
\int_{R}^{r}t^{1-N}dt\left( \int_{0}^{u\left( R\right) }f_{1}\left(
z,z\right) dz\right) ^{-1/2}+\overline{c}_{1}\overline{P}_{1\varepsilon }+%
\frac{1}{\varepsilon R^{\varepsilon }}\right) ,
\end{equation*}%
and 
\begin{equation*}
v\left( r\right) \geq a_{2}+\underline{Q}\left( r\right) .
\end{equation*}%
Therefore, if $\overline{P}_{1\varepsilon }\left( \infty \right) <\infty $
and $\underline{Q}\left( \infty \right) =\infty $ we get 
\begin{equation*}
\lim_{r\rightarrow \infty }u\left( r\right) <\infty \text{ and }%
\lim_{r\rightarrow \infty }v\left( r\right) =\infty \text{,}
\end{equation*}%
which yields the claim.

\textbf{2): }In this case, the idea is to mimic the proof of the Case 1. We
observe that%
\begin{equation}
u\left( r\right) \geq a_{1}+\underline{P}\left( r\right) \text{,}  \label{t2}
\end{equation}%
and%
\begin{equation}
v\left( r\right) \leq \mathcal{KO}_{f_{2}}^{-1}\left( \sqrt{C_{2}}%
\int_{R}^{r}t^{1-N}dt\left( \int_{0}^{v\left( R\right) }f_{2}\left(
z,z\right) dz\right) ^{-1/2}+\overline{c}_{2}\overline{P}_{2\varepsilon }+%
\frac{1}{\varepsilon R^{\varepsilon }}\right) .  \label{334}
\end{equation}%
The conclusion of the Theorem it follows by letting $r\rightarrow \infty $
in (\ref{t2}).

\subsection{\textbf{Proof of Theorem \protect\ref{th4} completed: }}

It follows from (\ref{ints}) and the conditions of the theorem that%
\begin{equation*}
\mathcal{KO}_{f_{1}}\left( u_{n}\left( r\right) \right) \leq \sqrt{2%
\overline{c}_{1}}\cdot \overline{P}_{1}\left( \infty \right) <\sqrt{2%
\overline{c}_{1}}\cdot \mathcal{KO}_{f_{1}}\left( \infty \right) <\infty ,
\end{equation*}%
and $v_{n}\left( r\right) \leq \mathcal{KO}_{f_{2}}^{-1}\left( \sqrt{2%
\overline{c}_{2}}\cdot \overline{P}_{2}\left( \infty \right) \right) <\infty 
$. On the other hand, since $\mathcal{KO}_{f_{1}}^{-1}$ is strictly
increasing on $\left[ 0,\mathcal{KO}_{f_{1}}\left( \infty \right) \right) $,
we find that $u_{n}\left( r\right) \leq \mathcal{KO}_{f_{1}}^{-1}\left( 
\sqrt{2\overline{c}_{1}}\cdot \overline{P}_{1}\left( \infty \right) \right)
<\infty $, and then the non-decreasing sequences $\left\{ u_{n}\right\}
_{n\geq 0}$ and $\left\{ v_{n}\right\} _{n\geq 0}$ are bounded above for all 
$r\geq 0$ and all $n$. We now use this observation to conclude $\left(
u_{n}\left( r\right) ,v_{n}\left( r\right) \right) \overset{n\rightarrow
\infty }{\rightarrow }\left( u\left( r\right) ,v\left( r\right) \right) $
and then the limit functions $u$ and $v$ are positive entire bounded radial
solutions of system (\ref{11}).\textbf{\ }

\subsection{\textbf{Proof of Theorem \protect\ref{th5} completed:}}

The proof for these cases is similar as the above and is therefore omitted.

\end{document}